\documentclass[11pt]{article}
\usepackage[auth-sc]{authblk}
\usepackage{amssymb}
\usepackage{amsthm}
\usepackage{graphicx}
\usepackage{latexsym}
\usepackage{amsmath,amscd}
\usepackage{bm}
\def\part#1{\frac{\partial\phantom{q}}{\partial#1}}

\newenvironment{rmk}{\begin{trivlist}\item[]{\bf Remark:} }
{\end{trivlist}}
\newenvironment{ex}{\begin{trivlist}\item[]{\bf Example:} }
{\end{trivlist}}
\newenvironment{rmks}{\begin{trivlist}\item[]{\bf Remarks:} }
{\end{trivlist}}
\newenvironment{exs}{\begin{trivlist}\item[]{\bf Examples:} }
{\end{trivlist}}
\newenvironment{prf}{\begin{trivlist}\item[]{\bf Proof:} }
{\hfill $\Box$ \end{trivlist}}
\newenvironment{lemprf}{\begin{trivlist}\item[]{\bf Proof:} }
 {\end{trivlist}}
\newtheorem{thm}{Theorem}

\newtheorem{prp}[thm]{Proposition}
\newtheorem{lemma}[thm]{Lemma}

\newcommand{\lie}[1]{\mathfrak{#1}}

\def\im{\mathop{\rm Im}\nolimits}

\newcommand{\R}{\mathbf{R}}
\newcommand{\C}{\mathbf{C}}
\newcommand{\K}{\mathbf{H}}
\newcommand{\Z}{\mathbf{Z}}

\newcommand{\CP}{{\mathbf C}{\rm P}}
\newcommand{\PP}{{\mathbf {\rm P}}}

\newcommand{\KP}{{\mathbf H}{\rm P}}
\title{ Manifolds with holonomy $U^*(2m)$}
 \author{Nigel Hitchin}

 \affil  {Mathematical Institute,
Radcliffe Observatory Quarter,
Woodstock Road,
Oxford, OX2 6GG}

\begin{document}
 \maketitle
 \thispagestyle{empty}

\let\oldthefootnote\thefootnote
\renewcommand{\thefootnote}{\fnsymbol{footnote}}
\footnotetext{hitchin@maths.ox.ac.uk}
\let\thefootnote\oldthefootnote
\noindent {\bf Abstract} We consider  the geometry determined by a torsion-free affine connection whose holonomy lies in the subgroup $U^*(2m)$, a real form of $GL(2m,\C)$,  otherwise denoted by $SL(m,\K)\cdot U(1)$.  We show in particular how examples may be generated from quaternionic K\"ahler or hyperk\"ahler manifolds with a circle action.
\vskip .25cm
\section{Introduction}
This article concerns the differential geometry of complex manifolds. The notion of a K\"ahler metric -- a Riemannian metric whose Levi-Civita connection is compatible with the complex structure -- has been of fundamental importance in the interaction between algebraic geometry and differential geometry. It can be considered as the geometry determined by a torsion-free affine connection whose holonomy lies in the subgroup $U(n)\subset GL(n,\C)$. Here we shall consider the geometry determined by an affine connection with holonomy in another real form, the group $U^*(2m)\subset GL(2m,\C)$. 

This structure has been little studied in the past. It made an appearance in the list of irreducible holonomy groups of Berger in 1955 \cite{Berg} but seems to have been overlooked (though not its complexification) in the complete  classification of Merkulov and Schwachh\"ofer \cite{Merk}. Only Joyce \cite{Joy1} and Pontecorvo \cite{Pont1} have touched on it. 

The author's interest derives from a geometric construction of Haydys \cite{Hay},  a generalization of the physicist's c-map construction, linking together other geometric structures based on the quaternions. It points to the existence of a large number of naturally occurring examples related one way or another to moduli spaces of interest to string theorists. However, we approach the subject from the point of view of a differential geometer. 

We first deal with the case of complex surfaces, where the structure takes the more familiar form of a K\"ahler metric of zero scalar curvature, about which there is considerable literature, which we summarize. Then we consider what properties as complex manifolds higher-dimensional compact examples  have, as a guide to the search for candidates. In particular, we show that the plurigenera vanish.

We then move to a construction of (non-compact) examples starting from a more familiar differential-geometric structure -- a quaternionic K\"ahler manifold of dimension $4m$. We show that on an open set of such a manifold endowed with a circle action there exists a natural circle-invariant connection with holonomy $U^*(2m)$. Following that we give a brief account of Haydys's result generating a quaternionic K\"ahler manifold with circle action from a hyperk\"ahler manifold with circle action. Since there are by now many examples of such hyperk\"ahler metrics, putting the two constructions together gives our suggested source of examples of $U^*(2m)$-manifolds.

Finally we treat the case of flat $\K^m$ with the action of right multiplication as a hyperk\"ahler manifold and describe the corresponding $U^*(2m)$-manifold. It turns out to be a product of $\CP^{m-1}$ with the quotient of a certain domain in $\C^{m+1}$ by an infinite cyclic group. This should provide a starting point for further examples in the future by the application of the quotient construction in \cite{Joy1}.

The work described here is based on the author's Santal\'o Lecture delivered in the Universidad Complutense on October 10th 2013. He would like to thank the organizers and ICMAT Madrid for support during the preparation of this paper.

\section{The holonomy group $U^*(2m)$}
In the classification of simple Lie groups the group $SU^*(2m)$ appears as a real form of $SL(2m,\C)$. It is defined in \cite{Helg} as the set of $2m\times 2m$ complex matrices of  determinant $1$  which commute with the antilinear map  
$$(z_1,\dots,z_m,w_1,\dots, w_m)\mapsto(\bar w_1,\dots, \bar w_m,-\bar z_1,\dots,-\bar z_m).$$ 
We define $U^*(2m)$ as the group obtained  by adjoining scalar multiplication by $e^{i\theta}$, so 
$$U^*(2m)=SU^*(2m)\times_{\pm 1}U(1)\stackrel{{\mathrm{def}}}=SU^*(2m)\cdot U(1)$$
and this is a real form of $GL(2m,\C)$. 

Other authors  have assigned a different  interpretation for the symbol $U^*(2m)$, but there is also another notation for $SU^*(2m)$, namely $SL(m,\K)$, the group of $m\times m$ quaternionic matrices with real determinant $1$.   Indeed, the antilinear map above is right multiplication by the quaternion $j$ on the quaternionic vector $w+jz$. The vector space $\C^{2m}$ may then be regarded as $\K^m$ with left multiplication by $SL(m,\K)$ and right multiplication by a unit complex number inside the quaternions.  This means that $U^*(2m)\cong SL(m,\K)\cdot U(1)$, a terminology  quite commonly used by differential geometers. However, for this paper we shall keep the notation $U^*(2m)$ by analogy with the compact real form $U(2m)$. 

We are interested in manifolds with a torsion-free connection whose holonomy lies in the real form $U^*(2m)$. This is one of the groups in Berger's original 1955 list \cite{Berg}. It has the property that for any reduction of the structure group of the tangent bundle to $U^*(2m)$, if there exists a compatible torsion-free connection, then it is unique \cite{Joy1},\cite{Sal1} (we shall also see this in the proof of Lemma \ref{11}). 

A special case is the group  $Sp(m)$, the  unitary quaternionic matrices, which is the maximal compact subgroup of $SU^*(2m)$ and these are called  {\it hyperk\"ahler manifolds}. Interestingly, the maximal compact subgroup $Sp(m)\cdot U(1)$ of $U^*(2m)$ is not an admissible holonomy group. Joyce in \cite{Joy1} calls holonomy $U^*(2m)$-manifolds {\it quaternionic complex} manifolds, but we want to downplay that quaternionic aspect here and focus on the complex one. We cannot avoid however remarking that the quaternionic point of view means that there is a twistor space interpretation of this geometry which gives it another degree of naturality.

This aspect derives from the study of manifolds $M^{4m}$ with torsion-free connections having holonomy in $GL(m,\K)\cdot\K^*$. They  are called {\it quaternionic manifolds} and Salamon \cite{Sal1} described their structure in terms of the complex structure of the twistor space, which is  a $2$-sphere bundle over $M$. Since $U^*(2m)=SL(m,\K)\cdot U(1)\subset GL(m,\K)\cdot\K^*$ our spaces are special cases. 

For a quaternionic manifold the twistor space is a complex manifold of dimension $2m+1$ with a real structure and a family of rational curves ({\it twistor lines}) with normal bundle isomorphic to ${\mathcal O}(1)^{2m}$ -- the direct sum of $2m$ copies of the line bundle ${\mathcal O}(1)$ of degree $1$ on $\CP^1$. The $2$-sphere fibres of the projection to $M$ are the twistor lines preserved by the real structure.  
The holomorphic structure of the normal bundle means that the anticanonical bundle $K^*=\Lambda^{2m+1}T$ of the twistor space restricts to ${\mathcal O}(2m+2)$ on each curve and what is required for a reduction to $U^*(2m)$ holonomy is a real section of a fractional power $(K^*)^{1/(m+1)}$ whose divisor  splits into two disjoint complex conjugate submanifolds. Each one of these intersects a twistor line in a single point and so provides a section of the fibration   and identifies $M$ with a complex manifold. The other component is the conjugate complex structure.

We shall need this viewpoint when we describe later the complex structure on certain concrete examples. Twistor theory, when it applies,  is a very useful coordinate-free method of telling us what is there but not so helpful in determining concrete formulae.

 The group $SL(1,\K)$ consists of simply the unit quaternions $Sp(1)$, and the left action is the action of $SU(2)$ on $\C^2$. Together with the right action of $U(1)$ this means that $U^*(2)=U(2)$. A manifold with this holonomy group is then just a K\"ahler manifold, but for consistency with the higher dimensional cases we make the convention that it is a K\"ahler manifold with {\it zero scalar curvature}. This is because the integrability of the twistor space in four dimensions is equivalent to the Weyl curvature being anti-self-dual and this condition for a K\"ahler metric is the vanishing of the scalar curvature.  Equivalently, as in \cite{Pont}, one can describe this via a holomorphic section of $(K^*)^{1/2}$ on the twistor space.    
 
\section{The case $m=1$}
Since two complex dimensions  is a special case, and also the only one which has been studied at all in depth, we discuss this next.

When a torsion-free connection has holonomy a subgroup of $GL(n,\C)$ then the almost complex structure is integrable and we have, by the Newlander-Nirenberg theorem, local complex coordinates $z_1,\dots, z_n$. A K\"ahler metric is Hermitian and so has the local form 
$$g=\sum_{\alpha,\beta} g_{\alpha\beta}dz_{\alpha}d\bar z_{\beta}$$
and the K\"ahler condition is that the $2$-form 
$$\omega=\frac{i}{2}\sum_{\alpha,\beta}  g_{\alpha\beta}dz_{\alpha} \wedge d\bar z_{\beta}$$
is closed. The Ricci form $\rho$ is the $2$-form defined locally by
$$\rho=-i\partial\bar\partial \log\det g_{\alpha\beta}$$
which is in fact coordinate-independent. Note that this local formula involves only the Hermitian metric $\det g_{\alpha\beta}$ on the anticanonical line bundle $K^*=\Lambda^nT$ and not the full Riemann curvature tensor. 

Writing 
$$\rho=-i\sum_{\alpha,\beta}  \rho_{\alpha\beta}dz_{\alpha} \wedge d\bar z_{\beta},$$
the scalar curvature is 
$$R=\sum_{\alpha,\beta}  g^{\alpha\beta}\rho_{\alpha\beta}.$$

We are concerned with the case of zero scalar curvature. Here is an example.
\begin{ex} Let $M=D \times \CP^1$ where the 2-sphere $\CP^1$ has its round metric of scalar curvature $+1$ and $D$ is the unit disc with the hyperbolic metric of constant curvature $-1$. In one dimension any metric is K\"ahler so the product is  K\"ahler and we have scalar curvature $1-1=0$. 
\end{ex}

This can be adapted to give compact examples. Let $\Sigma$ be a Riemann surface of genus $g>1$. By uniformization it has a hyperbolic metric induced from a discrete action of the fundamental group $\Gamma=\pi_1(\Sigma)$ on $D$. This is a  homomorphism from $\Gamma$  to $PSL(2,\R)$. Take now a representation of $\Gamma$ in $SO(3)$, the isometry group of $\CP^1$, and form the compact quotient
$$M=D\times_{\Gamma}\CP^1.$$
This is a $\CP^1$-bundle over $\Sigma$ with a scalar-flat K\"ahler metric: in fact from the theorem of Narasimhan and Seshadri it is identical with the projective bundle of a stable vector bundle. 

The $\CP^1$-fibres make $M$  into a ruled surface. This is one of the classes in  the Enriques-Kodaira classification of complex surfaces and we can appeal to the classification to restrict the candidates for compact surfaces admitting K\"ahler metrics of zero scalar curvature. The classification proceeds by considering the plurigenera of the surface: the dimension for large $N$  of the space of holomorphic sections of $K^N$. But we have:
\begin{prp} The plurigenera all vanish on a compact K\"ahler surface of zero scalar curvature unless the Ricci form vanishes identically.
\end{prp}
There are several ways of proving this but in the next section we shall give a proof which works for compact $U^*(2m)$-manifolds in all dimensions. A consequence of the proposition is that any such surface is either a K3 surface, an Enriques surface or a torus (the cases for which   the Ricci form vanishes) or 
is rational or ruled. 

This classification is up to birational equivalence, which means that candidates  could be obtained by blowing up points. A strong result in this direction is due to Kim, LeBrun and Pontecorvo \cite{KLP}, that any blow-up of a K\"ahler surface with $\rho\ne 0$ whose integral of the scalar curvature is non-negative has blow-ups which admit scalar-flat K\"ahler metrics. Note that the blow-ups are not arbitrary -- in fact blowing up points on a cubic curve in $\CP^2$ leads to a surface with a non-zero holomorphic section of $K^*$ and on a scalar-flat K\"ahler manifold these must vanish by the same argument as that used for the vanishing of the plurigenera. 

A full classification is not yet known and current interest in these metrics is driven by the conjectures of Donaldson, Tian and Yau relating the existence of a K\"ahler metric with constant scalar curvature metric, not necessarily zero, to a stability condition for projective varieties. 

\section{Compact $U^*(2m)$-manifolds}\label{compact}
The group $U^*(2m)$ for $m>1$ is non-compact, so in principle there is no underlying metric structure on a $U^*(2m)$-manifold $M^{4m}$. As a subgroup of $GL(2m,\C)$, the torsion-free connection defines the structure of a  $2m$-dimensional complex manifold. The subgroup $SU^*(2m)$ lies in $SL(2m,\C)$, which acts trivially on the canonical bundle $K$ and so the $U(1)$ factor in $U^*(2m)$ defines a Hermitian structure on $K$, with a Ricci form but no metric to contract and form a scalar curvature. Nevertheless the Hermitian structure  provides a non-vanishing real section of $K\bar K$, a volume form. And, as with all affine connections, we have the notion of geodesic.

In higher dimensions it is far more difficult as yet to find compact examples of manifolds with holonomy $U^*(2m)$. We can look at subgroups,  and compact hyperk\"ahler manifolds with holonomy $Sp(m)$ do exist. There are relatively few of these however -- Hilbert schemes of K3 surfaces or tori and two exceptional cases due to O'Grady \cite{OG}. The subgroup $SU^*(2m)$ corresponds to the case of vanishing Ricci form and there are a few examples here. Swann \cite{swann2} has a construction which gives simply-connected examples  based on quotients of torus bundles over K3 surfaces and Barberis et al. \cite{Barb} produce compact nilmanifolds of this type. These manifolds are far from algebraic, in fact Verbitsky has shown \cite{Verb} if $M$ admits any K\"ahler metric then in this case it admits a hyperk\"ahler metric.  We thus have the analogues of the K3 surface and the torus (in the sense that the Ricci form vanishes) but nothing like the  ruled surfaces. Nevertheless, in \cite{Pont1} Pontecorvo shows the existence of compact $12$-dimensional examples where the Ricci form does not vanish, though being symmetric spaces their holonomy is not the full group $U^*(2m)$. 

Despite the paucity of examples it still makes sense  to look for  complex manifolds which might be candidates. 
There is a rough classification for higher dimensional complex manifolds as initiated by Ueno \cite{Ueno} and the plurigenera play a role here, so the following proposition helps to narrow the field:

\begin{prp} The plurigenera all vanish on a compact $U^*(2m)$-manifold  unless the Ricci form vanishes identically.
\end{prp} 

\begin{prf} We first show that the Ricci form of a $U^*(2m)$-manifold has a specific algebraic form. We then   use  a vanishing theorem to prove the result. 

To proceed  we revert to the $SL(m,\K)\cdot U(1)$ description: choosing a local trivialization of the $U(1)$-bundle,  we have  locally defined almost complex structures $J$ and $K$.
\begin{lemma}\label{11} The Ricci form $\rho$ of a $U^*(2m)$-manifold is of type $(1,1)$ with respect to all complex structures $I,J,K$.
\end{lemma}
\begin{lemprf} The determinant of $U^*(2m)$ induces a $U(1)$ connection on $K^*=\Lambda^{2m}T^{1,0}$ compatible with the  holomorphic structure and is hence the canonical Chern connection whose curvature, the Ricci form $\rho$, is of type $(1,1)$ with respect to the complex structure $I$.

For $J$ and $K$ we need to refer to Salamon \cite{Sal1} for a description of the curvature tensor of a manifold with a torsion-free $GL(m,\K)\cdot \K^*$-connection. In this case there is an induced $SO(3)$-connection  on the bundle $Q$ of imaginary quaternions acting on the tangent bundle.  For  holonomy $U^*(2m)$, $I$ is a covariant   constant section of $Q$. Salamon's result is that the curvature tensor is of the form

\begin{equation}
\sum_i\partial (v_i\otimes t^i)+     R_U
\label{curv}
\end{equation}

where $R_U$ lies in an irreducible representation  of $SL(m,\K)$, $v_i$ lies in the first prolongation $\lie{g}^{(1)}$ of the Lie algebra $\lie{g}$ of $GL(m,\K)\cdot \K^*$ and $t_i\in T^*$. 
The bundle $\lie{g}^{(1)}\subset \lie{g}\otimes T^*$ is in this case isomorphic to $T^*$, and the homomorphism $\partial$  is the restriction of the natural map $\lie{g}\otimes T^*\otimes T^*\mapsto \lie{g}\otimes \Lambda^2T^*$.

More concretely   $\alpha\in T^*$   acts non-trivially on the tensors corresponding to three different representations of  $GL(m,\K)\cdot \K^*$ on its Lie algebra. Evaluated on a tangent vector $Z$, the first takes values in the scalars in $\lie{gl}(m,\K)$ and is  $Z\mapsto \alpha(Z) 1$;  the second, in $\lie{gl}(m,\K)$ is 
$$Z\mapsto - Z\otimes \alpha+IZ\otimes I\alpha+JZ\otimes J\alpha+KZ\otimes K\alpha $$
and the  third in $ \K$  is $Z\mapsto \alpha(IZ) I+\alpha(JZ) J+\alpha(KZ) K$. 

\begin{rmk} The first prolongation measures the choice in finding a torsion-free connection and in this case the above description of the action shows that any two  torsion-free connections $\nabla,\tilde\nabla$ preserving a $GL(m,\K)\cdot \K^*$-structure are related by a $1$-form $\alpha$ as follows:  
\begin{eqnarray}
\tilde\nabla_ZY&=&\nabla_ZY+ \alpha(Z)Y+\alpha(Y)Z- \alpha(IY)IZ-\alpha(IZ)IY- \nonumber\\
&-&\alpha(JY)JZ-\alpha(JZ)JY-\alpha(KY)KZ-\alpha(KZ)KY. \label{form}
\end{eqnarray}
 The action on the real line bundle $\Lambda^{4m}T$ is $(4m+2)\alpha$. This shows, as mentioned earlier, that a $U^*(2m)$-connection, when it exists, is unique, since it preserves a volume form and there is no more freedom to choose $\alpha$.
\end{rmk}

Now the $R_U$ term in the curvature acts trivially on the $SO(3)$ bundle $Q$ as do the first and second terms, hence we may assume that
$$v_i=I\otimes I\alpha_i+J\otimes J\alpha_i+K\otimes K\alpha_i.$$
It follows that  the curvature of $Q$ is 
$$\sum_iI\otimes I\alpha_i\wedge v_i+J\otimes J\alpha_i\wedge v_i+K\otimes K\alpha_i\wedge v_i$$
where $I,J,K$ act by $q\mapsto[I,q]$ etc. on the imaginary quaternion $q$. 

We know however that a $U^*(2m)$-connection preserves a global complex structure $I$ and hence the curvature commutes with $I$ and so 
\begin{equation}
\sum_i J\alpha_i\wedge v_i=0=\sum_iK\alpha_i\wedge v_i
\label{JK}
\end{equation}
which gives  the Ricci form  
$$\rho=\sum_i I\alpha_i\wedge v_i.$$

A $2$-form is of type $(1,1)$ with respect to $J$ if it is annihilated by the Lie algebra action $\alpha\wedge \beta\mapsto J\alpha\wedge \beta+\alpha\wedge J\beta$. Apply this action to $\rho$ as above and we get 
$$\sum_i -K\alpha_i\wedge v_i+\sum_i I\alpha_i\wedge Jv_i=\sum_i I\alpha_i\wedge Jv_i$$
from (\ref{JK}).
But applying $J\wedge J$ to the second equation of (\ref{JK}) we get 
$$\sum_i I\alpha_i\wedge Jv_i=0$$
and hence $\rho$ is of type $(1,1)$ with respect to $J$. Similarly with respect to $K$. \qed 

\begin{rmk} In Salamon's notation $\rho$ lies in the subbundle $S^2E$ of $\Lambda^2T^*$. It means that the connection is defined in twistor terms by a holomorphic line bundle on twistor space. In fact it is the line bundle corresponding to the difference of the disjoint divisors on which the holomorphic section of $(K^*)^{1/(m+1)}$ defining the $U^*(2m)$-structure vanishes.
\end{rmk}
\end{lemprf}
 To continue with the proof, choose a reduction of the structure group to its maximal compact subgroup $Sp(m)\cdot U(1)$. The connection does not reduce but we do have a Hermitian metric since $Sp(m)\cdot U(1)\subset SU(2m)\cdot U(1)=U(2m)$. To each local almost complex structure $J,K$ there are then  Hermitian forms $\omega_2,\omega_3$ and furthermore $\omega_3+i\omega_1$ is of type $(2,0)$ with respect to $J$. This means that $\omega_1$, the Hermitian form for $I$, is of type $(2,0)+(0,2)$ with respect to $J$. But from the Lemma 3, $\rho$ is of type $(1,1)$ with respect to $J$ and hence orthogonal to $\omega_1$, which means that the connection on the anticanonical bundle  is a Hermitian Yang-Mills connection. 
 
 We may now apply a standard vanishing theorem for such connections \cite{K} to deduce that any holomorphic section of $K^N$ must be covariant constant. Hence unless  $\rho=0$ the plurigenera vanish. Note that this applies to both positive and negative powers of the canonical bundle $K$ of $M$.
\end{prf}

\section{From quaternionic K\"ahler to holonomy $U^*(2m)$} \label{quatu}
It seems difficult to find examples of $U^*(2m)$-manifolds {\it ab initio} but we can in fact find some from more familiar differential-geometric territory. Here we shall show that a quaternionic K\"ahler manifold with a circle action, or more generally a Killing vector field,  generates such a structure on an open set. The four-dimensional version is the known conformal equivalence of self-dual Einstein and scalar-flat K\"ahler with circle symmetry \cite{Tod1}.

Let $M^{4m}$ be a quaternionic K\"ahler manifold with a circle action preserving the structure. By definition the holonomy group of a quaternionic K\"ahler manifold is $Sp(m)\cdot Sp(1)\subset GL(m,\K)\cdot\K^*$ and this is the Levi-Civita connection.  The rank $3$ bundle associated to the adjoint representation of $Sp(1)$ is the bundle $Q$ of imaginary quaternions we considered in the previous  section. With  the given metric we may also consider the skew adjoint transformations $I,J,K$ as $2$-forms and then we have a local oriented orthonormal basis $\omega_1,\omega_2,\omega_3$ of $Q$ consisting of local hermitian forms for $I,J,K$.

If $X$ is the vector field generated by the action then there is (see \cite{GL}) a  section of $Q$, the {\it moment section} ${\bm\mu}=\mu_1\omega_1+\mu_2\omega_2+\mu_3\omega_3$,  which satisfies 
\begin{equation}
\nabla {\bm\mu}=\sum_i i_X\omega_i\otimes \omega_i.
\label{mom}
\end{equation}
Given $X$, the most direct way to determine ${\bm \mu}$ is to take the exterior derivative of  the dual 1-form $X^{\flat}$: its   component in the subbundle  spanned by $\omega_1,\omega_2,\omega_3$ is a multiple of   ${\bm \mu}$. In four dimensions this is the self-dual component of the $2$-form 
$dX^{\flat}$.

On the complement  $M_0$ of its zero set,  ${\bm\mu}$  picks out a distinguished almost complex structure $I$.  It was shown by  Battaglia \cite{Batt} that   $I$ is integrable. We prove the stronger result: 

\begin{prp} There is a unique torsion-free affine connection with holonomy $U^*(2m)$ on $M_0$ defining the same quaternionic structure as its Levi-Civita connection. 
\end{prp}

\begin{prf}
 Let $\nabla$ be the  Levi-Civita connection  with holonomy in $Sp(m)\cdot Sp(1)$. Then, as in Equation (\ref{form})   above, we need to find a connection $\tilde\nabla$ of the form 
$$\tilde\nabla_ZY=\nabla_ZY+ \alpha(Z)Y+\alpha(Y)Z- \alpha(IY)IZ-\alpha(IZ)IY+\dots$$
which preserves $I$: to   determine an appropriate $1$-form $\alpha$.  
 
  Choose  a local basis such that ${\bm\mu}=\mu_1\omega_1$ then since  $\nabla \omega_1=\theta_2\otimes \omega_3-\theta_3\otimes \omega_2$ etc. for  forms $\theta_i$ which are components of the connection in this basis, the equation (\ref{mom}) for ${\bm\mu}=\mu_1\omega_1$ reads
\begin{equation}
d\mu_1=i_X\omega_1,\qquad
\mu_1\theta_2=i_X\omega_3,\qquad
\mu_1\theta_3=-i_X\omega_2.
\label{mom1}
\end{equation}
 
 In this basis $\nabla I=\theta_2\otimes K-\theta_3\otimes J$ and, as above,  $\alpha$ acts on $I$ as $J\alpha\otimes [J,I]+K\alpha\otimes [K,I]=-2J\alpha\otimes K+2K\alpha\otimes J$ so  $$\tilde\nabla I = \theta_2\otimes K-\theta_3\otimes J+2J\alpha\otimes K-2K\alpha\otimes J$$
  and to make this vanish we require    $\theta_2=-2J\alpha$ and $\theta_3=-2K\alpha$. 
 Note that these two equations are consistent,  for from (\ref{mom1}) we have $i_X(\omega_2+i\omega_3)=i\mu_1(\theta_2+i\theta_3)$ which is of type $(1,0)$ hence  $\theta_3=I\theta_2$ and so $\alpha=J\theta_2/2=K\theta_3/2$. 
 
 It follows that 
$$2\alpha(Y)=(J\theta_2)(Y)=\mu_1^{-1}\omega_3(X,JY)=\mu_1^{-1}g(KX,JY)=-\mu_1^{-1}\omega_1(X,Y)$$
and from (\ref{mom1}) this gives
$\alpha = -d\log \mu_1/2$ as the required $1$-form. We have already observed uniqueness from the preserved volume condition.

The holonomy  therefore lies in $U^*(2m)$.
\end{prf}

\begin{rmks}

\noindent 1. The Riemannian volume form $\nu$ is preserved by $\nabla$ but acted on by the real trace of the extra term in  (\ref{form}). This is $(4m+4)\alpha=(2m+2)d\log \mu_1$. Hence 
$$\tilde\nabla \nu=-(2m+2)(d\log \mu_1)\otimes \nu$$
and $\mu_1^{-(2m+2)}\nu$ is an invariant volume form. 

\noindent 2. When $m=1$, the formula (\ref{form}) for the torsion-free connection becomes 
$$\tilde\nabla_ZY=\nabla_ZY+ \alpha(Z)Y+\alpha(Y)Z-g(Y,Z)\alpha^{\sharp}$$
and in general this is valid if $Y,Z$ lie in the same one-dimensional  quaternionic subspace of the tangent space. The formula for $m=1$ and this choice of $\alpha$ is precisely the Levi-Civita connection for the metric  ${\mu_1}^{-2}g$. The interpretation of quaternionic K\"ahler in $4$-dimensions is of an anti-self-dual Einstein metric so here we have one which is conformally equivalent to a scalar-flat K\"ahler metric. This is described explicitly in  \cite{Dun}. 

\noindent 3. This proposition is basically the same as that of Joyce in \cite{Joy1} in the sense that our moment section is his twistor function in the rather more general setting of quaternionic manifolds. He proves that given a twistor function one has a volume form and hence a distinguished connection in the quaternionic family, then a calculation shows that it preserves a complex structure $I$. We, on the other hand,  start with a quaternionic K\"ahler manifold   and the above formula gives the $U^*(2m)$-connection explicitly in terms of the Levi-Civita connection. 
\end{rmks}
\begin{exs} 

\noindent 1. Let the quaternionic K\"ahler manifold  be the quaternionic projective line $\KP^1=S^4$. This is conformally flat by stereographic projection and writing $\R^4=\R^2\times \R^2$  
we have  
$$g_{S^4}=\frac{1}{(1+\rho^2+\sigma^2)^2}(d\rho^2+\rho^2d\varphi^2+d\sigma^2+\sigma^2d\theta^2)$$
Consider the rotation action on the right hand $\R^2$-factor. This is  generated by the vector field $X=\partial/\partial \theta$. It is an  isometric action and $X^{\flat}=\sigma^2 d\theta/(1+\rho^2+\sigma^2)^2$ which is $u^2d\theta$ for $u=\sigma/(1+\rho^2+\sigma^2)$. Changing coordinates to $(u,v,\varphi, \theta)$ where $v=(\rho^2+\sigma^2-1)/\rho$ gives, using the relation $(1+\rho^2-\sigma^2)^2+4\sigma^2\rho^2=(1+\rho^2+\sigma^2)^2-4\sigma^2=(\rho^2+\sigma^2-1)^2+4\rho^2$, the following expression:
$$g_{S^4}=\frac{1}{1-4u^2}du^2+ u^2d\theta^2+\frac{1-4u^2}{(v^2+4)^2}dv^2+\frac{1-4u^2}{v^2+4}d\varphi^2$$
Since  $dX^{\flat}=2udu\wedge d\theta$ the self-dual component of $dX^{\flat}$ is 
$$(dX^{\flat})^+= udu\wedge d\theta+\frac{(1-4u^2)^{3/2}}{(v^2+4)^{3/2}}dv\wedge d\varphi.$$      This is the moment section ${\bm\mu}$.    
Hence 
$$\mu_1^2\omega_1^2=((dX^{\flat})^+)^2=2u\frac{(1-4u^2)^{3/2}}{(v^2+4)^{3/2}}du\wedge d\theta\wedge dv\wedge d\varphi=(1-4u^2)vol_{S^4}.$$
and  the metric 
$$g=\frac{1}{(1-4u^2)^2}du^2+ \frac{1}{(1-4u^2)}u^2d\theta^2+\frac{1}{(v^2+4)^2}dv^2+\frac{1}{v^2+4}d\varphi^2$$
is a scalar-flat K\"ahler metric. In fact it is the product of a hyperbolic plane of curvature $-4$ and a sphere of curvature $+4$ as can be seen by substituting $u=(\tanh 2x)/2$ and $v=2 \tan y$. This is well-defined for $4u^2\ne 1$, i.e. the complement of the circle $\rho=0,\sigma=1$ in $S^4$.

\noindent 2. The above is an explicit calculation but we can extend it to higher dimensions by considering quaternionic projective space $\KP^m=Sp(m+1)/Sp(m)\cdot Sp(1)$ which is a quaternionic K\"ahler manifold with its symmetric space metric. There is a circle action defined by left multiplication by $e^{i\theta}$ on $\K^{m+1}$. The complex structure on the associated $U^*(2m)$-manifold is best seen via the twistor space which for quaternionic projective space is just complex projective space $\CP^{2m+1}$. The circle action decomposes $\K^{m+1}$ as a right complex module as $V\oplus V^*$ and allows us to write the twistor space  as $\PP(V\oplus V^*)$ and then the section of $(K^*)^{1/(m+1)}\cong {\mathcal O}(2)$ is the natural quadratic pairing $\langle z,w\rangle$ between $z\in V$ and $w\in V^*$. Removing the subset $\Vert z\Vert^2-\Vert w\Vert^2=0$ gives two disconnected components, so as a complex manifold our $U^*(2m)$-manifold is 
$$\{(z,w)\in \PP(V\oplus V^*): \langle z,w\rangle=0; \Vert z\Vert^2-\Vert w\Vert^2>0\}$$
an open set in a projective quadric. Recall that a $2$-dimensional quadric is $\CP^1\times \CP^1$ and removing $S^1\times \CP^1$ gives the scalar-flat K\"ahler metric on the component $D\times \CP^1$.

This example is a non-Riemannian symmetric space $SO^*(2m+2)/SO^*(2m)\times SO(2)$, a noncompact form of the Grassmannian $SO(2m+2)/SO(2m)\times SO(2)$ which is well known to be the projective complex quadric.    
The group $SO^*(2m)$  is the subgroup of $GL(2m,\C)$ which preserves a complex inner
product $(u,v)$ and commutes with an antilinear automorphism $J$ with $J^2=-1$   for which
$(u, Jv)$ is a hermitian form. If ($u, Ju) > 0$ then ($Ju, J^2u) = -(Ju, u) < 0$  and so the
form has hermitian signature $(m,  m)$. Because $J$ defines a quaternionic structure the alternative notation is $SO(m,\K)$.

This is the basis for  the example of Pontecorvo  \cite{Pont1}: the case $m=3$ has a compact quotient  which is his compact 12-dimensional example. These are locally symmetric, and in fact, as pointed out in \cite{Pont1} they have an  indefinite K\"ahler structure: a torsion-free connection with holonomy $U(m,m)$. Note that for $m=1$ the manifold is $D\times \CP^1$ and the $U^*(2)$-connection is a product. Since $U(1)\times U(1)=U(2)\cap U(1,1)$  the connection also preserves a $U(1,1)$-structure. 
\end{exs}

Every compact simple group $G$ has an associated quaternionic K\"ahler manifold, its  Wolf space $G/K\cdot Sp(1)$, and consequently many circle actions so there are plenty of noncompact examples to be had by this construction. It is also possible to adapt the construction to indefinite metrics --  $Sp(p,q)\cdot Sp(1)$-holonomy -- since $U^*(2m)$ does not respect any particular signature. There is also the quaternionic K\"ahler quotient \cite{GL} which, applied to  $\KP^n$, usually produces orbifold singularities. It is possible, nevertheless, that these lie on the zero set of the moment section and are not inherited by the associated $U^*(2m)$-manifold.

However, to make  effective progress we should be able to generate more quaternionic K\"ahler metrics. 
\section{From hyperk\"ahler to holonomy $U^*(2m)$}
If we lack naturally occurring quaternionic K\"ahler manifolds, the same is not true of their cousins the hyperk\"ahler manifolds. Despite the relatively few compact examples, as mentioned in Section \ref{compact} we know many examples of noncompact ones. These appear on certain moduli spaces of gauge-theoretic equations and also through the hyperk\"ahler quotient construction. 

A   result of Haydys  \cite{Hay}  shows that a hyperk\"ahler manifold with a certain type of circle action generates a quaternionic K\"ahler manifold with circle action and hence, by the results of  Section \ref{quatu}, a manifold with holonomy $U^*(2m)$. This construction is a generalization of the c-map in physics where the initial hyperk\"ahler manifold is associated to a gauge theory in four dimensions with $N=2$ supersymmetry.

To see how to generate examples of hyperk\"ahler manifolds, we first  briefly describe the hyperk\"ahler quotient \cite{HKLR}. Recall that a hyperk\"ahler manifold $M$ has globally-defined complex structures $I,J,K$, behaving like quaternions, and a metric $g$ which provides corresponding K\"ahler forms $\omega_1,\omega_2,\omega_3$. Suppose we have a free action of a Lie group $G$ which preserves all three K\"ahler forms, and has well-defined  equivariant moment maps $\nu_1,\nu_2,\nu_3$ taking values in the dual $\lie{g}^*$ of the Lie algebra, or equivalently a vector-valued moment map ${\bm \nu}:M\rightarrow \lie{g}^*\otimes \R^3$. Then the result is that the quotient metric on ${\bm \nu}^{-1}(0)/G$ is hyperk\"ahler. 

\begin{ex}  The simplest example is to take flat space $\K^{m+1}=\C^{m+1}\oplus j\C^{m+1}$ and $G=U(1)$ where the action is left multiplication by $e^{i\theta}$. In the complex structure $I$, with $V=\C^{m+1}$, this is $V\oplus V^*$ with the action $(z,w)\mapsto (e^{i\theta}z,e^{-i\theta}w)$. The basic moment map is then $\nu_1=(\Vert z\Vert^2-\Vert w\Vert^2)/2$ and $\nu_2+i\nu_3=\langle v,w\rangle$, but we can add constants and keep equivariance since $G$ is abelian. In particular if  $\nu_1=(\Vert z\Vert^2-\Vert w\Vert^2)/2-1$ then on the zero set $z$ is non-vanishing and defines a map to the complex projective space $\PP(V)$. It is not hard to see then that the  condition $\nu_2+i\nu_3=\langle z,w\rangle=0$ leads to the quotient being the cotangent bundle ${\mathrm T}^*\PP(V)$. This is known as the {\it Calabi metric} \cite{Cal}.
\end{ex} 

To construct a quaternionic K\"ahler manifold we need a hyperk\"ahler manifold $M$ with another type of circle action: one which fixes $\omega_1$ but acts on $\omega_2+i\omega_3$ by multiplication by $e^{i\theta}$, for example $(z,w)\mapsto (z, e^{i\theta}w)$ in flat space. If $\mu$ is a moment map for the circle action using the form $\omega_1$ then one observes that $F=\omega_1+dd^c_1\mu$ is of type $(1,1)$ with respect to all complex structures \cite{Hay}. If the cohomology class of $F/2\pi$ is integral then $F$ is the curvature of a principal $U(1)$-bundle with connection over $M$. Given a lift of the circle action to $P$, Haydys shows that there is a natural quaternionic K\"ahler metric (which may have indefinite signature) on the quotient $\hat M$ of $P$ by the lifted action. The principal bundle $U(1)$-action now becomes a geometric action on $\hat M$ and from here we can construct on an open set a manifold of holonomy $U^*(2m)$ as in Section \ref{quatu}.  A description of this construction focused on the twistor approach is given in \cite{Hit}. 

The process is reversible and is called the hyperk\"ahler/quaternionic K\"ahler correspondence. It will be useful for future reference to describe now the reverse route. This uses  the  Swann bundle associated to any quaternionic K\"ahler manifold  \cite{swann}.  On a quaternionic K\"ahler manifold $\hat M$ we have the bundle of imaginary quaternions $Q$. Let $\hat P$ be its principal $SO(3)$ frame bundle. Let $E_1,E_2,E_3$ be the vector fields on $\hat P$ generated by the standard basis of $\lie{so}(3)$, then the connection form has components given by  $1$-forms $a_1,a_2,a_3$ on $\hat P$ where $i_{E_i}a_j=\delta_{ij}$.  The curvature of this connection  is given by
 $$da_1+a_2\wedge a_3=c\,\omega_1, \quad da_2+a_3\wedge a_1=c\,\omega_2,\quad da_3+a_1\wedge a_2=c\,\omega_3$$
 where $c$ is the constant scalar curvature of $\hat M$.

The $(4m+4)$-manifold $\hat P\times \R^+$ is  known as the Swann bundle of $\hat M$. If  the scalar curvature of $\hat M$ is positive,  then the triple of closed $2$-forms  $\varphi_i=d(ta_i)$ (where $t$ is the $\R^+$-coordinate) are the K\"ahler forms for a hyperk\"ahler metric. If the scalar curvature is negative then  the holonomy lies in $Sp(m,1)$ rather than $Sp(m+1)$. 
 
 Suppose now we have a circle action on $\hat M$ generated  by a vector field $X$ with a corresponding moment section ${\bm \mu}$ of $Q$. Then we can lift it naturally to an action on the Swann bundle as follows.   Let $\bar X$ denote the horizontal lift of  $X$ on $\hat M$ to $\hat P$, then $i_{\bar X}a_i=0$. Define a new lift by $Y=\bar X-c\sum_i\mu_iE_i$. Then expanding $\varphi_1$ as 
$$\varphi_1=d(ta_1)=dt\wedge a_1+ta_2\wedge a_3+tc\, \omega_1 $$
we obtain 
 $$i_Y\varphi_1=  tc(d\mu_1+\mu_2a_3-\mu_3 a_2)-c(-\mu_1dt+t\mu_2a_3-t\mu_3a_2)=c\,d(\mu_1t)$$
using $d\mu_1+\mu_2a_3-\mu_3 a_2=i_{\bar X}\omega_1$ from (\ref{mom}), and similarly for $\varphi_2,\varphi_3$.
  Then ${\mathcal L}_Y\varphi_i=0$ and the action preserves the three K\"ahler forms. Moreover  $c(\mu_1t,\mu_2t,\mu_3t)$ is a   hyperk\"ahler moment map.  
  
  We take a hyperk\"ahler quotient  of the Swann bundle by this lifted circle action to obtain the hyperk\"ahler manifold $M$. If we set ${\bm \nu}=c(\mu_1t,\mu_2t,\mu_3t)=(c,0,0)$ then ${\bm \nu}$ is preserved by a circle subgroup  $SO(2)\subset SO(3)$, with generating vector field $E_1$, acting on $\hat P$ and this descends to $M$ as the required circle action.

 \begin{rmk} We can calculate the  moment map $\mu$ for the action on $M$ by $d\mu=i_{E_1}\varphi_1=-dt$ and so, since $\mu_1t=1$ on  ${\bm \nu}^{-1}(0)$, 
 $$\mu=-t=-\frac{1}{\mu_1}.$$
Hence, approaching the zero set of ${\bm \mu}$ on the quaternionic K\"ahler manifold $\hat M$ corresponds to $\mu$ on $M$ going to infinity. 
\end{rmk}

\begin{ex} Start with the hyperk\"ahler manifold $M=T^*\CP^m$ and its Calabi metric, together with  the circle action given by scalar multiplication by $e^{i\theta}$ on the fibres. Then we claim that this generates by the correspondence and the construction of Section \ref{quatu} the symmetric space $SO^*(2m+2)/SO^*(2m)\times SO(2)$ as a $U^*(2m)$-manifold. Note that although the zero section of the cotangent bundle is fixed by the circle action, in fact the natural lifted action on the $U(1)$-bundle is non-trivial there so we do not need to remove the zero section to implement the construction. 

To prove that this is the correspondence, we note that the quaternionic K\"ahler manifold $\KP^m$ with the action of the diagonal circle in $Sp(m+1)$ gave us in Section  \ref{quatu} the $U^*(2m)$-manifold $SO^*(2m+2)/SO^*(2m)\times SO(2)$. On the other hand the Swann bundle of $\KP^m$ is actually the flat hyperk\"ahler manifold $\K^{m+1}\backslash \{0\}$. To get the hyperk\"ahler manifold corresponding to $\KP^m$ we take the hyperk\"ahler quotient of the lifted action at a nonzero value of the moment map. But the diagonal circle group acts on $\K^{m+1}\backslash \{0\}$ exactly as in the example  above which gave the Calabi metric.
\end{ex}
Since most known hyperk\"ahler metrics arise as a quotient of flat space we next investigate which $U^*(2m)$-manifold corresponds to $\K^m$ with its standard circle action.

\section{Flat space}
\subsection{Four dimensions}\label{41}
In four dimensions $U^*(2)$-geometry is a metric geometry and hence it is easier to write down just what a $U^*(2)$-structure is. In the next section we shall deal with higher dimensions by focusing only on the complex structure. 

Flat space as a hyperk\"ahler manifold is defined by  $M=\C^2$ (complex structure $I$) and 
 $$\omega_1=\frac{i}{2}(dz \wedge d\bar z+dw \wedge d\bar w), \qquad \omega_2+i\omega_3=dz\wedge dw.$$
Then  $(z,w)\mapsto (z,e^{i\theta}w)$ takes  $\omega_2+i\omega_3$ to $e^{i\theta}(\omega_2+i\omega_3)$ and is the circle action we consider. This gives the vector field $$X=iw\frac{\partial}{\partial w}-i\bar w\frac{\partial}{\partial \bar w}$$
 and the moment map 
 $\mu=-\vert w\vert^2/2$, 
so $dd_1^c\mu=-idw\wedge d\bar w$ and  $$F=\omega_1+dd_1^c\mu=\frac{i}{2}(dz\wedge d\bar z-dw \wedge d\bar w).$$
This is the curvature of the trivial holomorphic line bundle with the connection compatible with the Hermitian metric $\exp(\vert z\vert^2-\vert w\vert^2)/2$. We take the lifted  action to be trivial, which means that we must remove the fixed point set $w=0$ to implement the correspondence. It follows that the quaternionic K\"ahler manifold -- an $S^1$-quotient of $\C\times \C^*\times S^1$ -- has infinite cyclic fundamental group.

An explicit expression for this metric is given in \cite{ACDM} as a 4-manifold with an isometric action of a 3-dimensional   group  $H$. The abelian group of translations $(z,w)\mapsto (z+c,w)$ acts as hyperk\"ahler isometries of $\C^2$, commuting with the circle action, but the constant curvature form $F$ means that it lifts to act on the principal $U(1)$-bundle via a central extension --  a Heisenberg group. 

This becomes a geometric action on the quaternionic K\"ahler  side. 
Let  $\sigma_i$ be invariant 1-forms on $H$ dual to the standard generators, then  

\begin{equation}
d\sigma_1=2\sigma_2\wedge \sigma_3\qquad d\sigma_2=d\sigma_3=0
\label{halg}
\end{equation} 

and the quaternionic K\"ahler metric is 
$$g=\frac{1}{4\rho^2}\left(d\rho^2 + \sigma_1^2+2\rho(\sigma_2^2+\sigma_3^2)\right)$$
for $\rho>0$.

This metric is well-known -- its universal covering is the Bergman metric on the unit ball in $\C^2$ and the Heisenberg orbits are horospheres through a point on the boundary.  Alternatively it  is the Bergman metric on the biholomorphically equivalent  Siegel domain $\{(v_1,v_2)\in \C^2: \im v_2>\vert v_1\vert^2\}$. In this model  $\rho= \im v_2-\vert v_1\vert^2$  and $\hat M$ is the quotient by the $\Z$-action $(v_1,v_2)\mapsto (v_1,v_2+2\pi n)$.

The invariant $1$-form $\sigma_1$ is dual to the centre of $H$ which was the principal bundle action on the $U(1)$-bundle over $M$ which now becomes the geometric action on $\hat M$, so  $X^{\flat}=\sigma_1/4\rho^2$ and 
$$dX^{\flat}=-\frac{1}{2\rho^3}d\rho\wedge \sigma_1+ \frac{1}{2\rho^2}\sigma_2\wedge \sigma_3. $$
Since  $\ast d\rho\wedge \sigma_1=2\rho\sigma_2\wedge \sigma_3$ we have 
\begin{equation}
(dX^{\flat})^+= -\frac{1}{8\rho^3}d\rho\wedge \sigma_1- \frac{1}{4\rho^2}\sigma_2\wedge \sigma_3. 
\label{quatI}
\end{equation}

so that 
$$((dX^{\flat})^+)^2=\frac{1}{16\rho^5}d\rho\wedge\sigma_1\wedge\sigma_2\wedge\sigma_3=\frac{1}{\rho^2}vol_g. $$
Thus, following Section \ref{quatu},
\begin{equation}
\tilde g=\frac{1}{4}\left(d\rho^2 + \sigma_1^2+2\rho(\sigma_2^2+\sigma_3^2)\right)
\label{scal}
\end{equation}
is a scalar-flat K\"ahler metric. 

From (\ref{quatI}) and the definition of the complex structure $I$ in Section \ref{quatu} the $(1,0)$-forms are spanned by  $d\rho+i\sigma_1, \sigma_2+i\sigma_3$ so $\bar\partial \rho^2=2\rho (d\rho)^{0,1}=\rho(d\rho-i\sigma_1)$, hence
$$\partial\bar\partial \rho^2=-id(\rho\sigma_1)=-id\rho\wedge
\sigma_1-2i\rho\sigma_2\wedge\sigma_3$$ 
and $\rho^2/2$ is a K\"ahler potential. 

One can check that $I$ is {\it not} the complex structure of the Siegel domain, but the metric is simple enough that we can calculate it.   Still working on the universal cover we have  $d(\sigma_2+i\sigma_3)=0$ from  Equation (\ref{halg}) so $\sigma_2+i\sigma_3= dw_1$ for some complex function $w_1$. However  $\sigma_2+i\sigma_3$ is of type $(1,0)$ so $w_1$ is holomorphic. 
Now $d\sigma_1=2\sigma_2\wedge\sigma_3=idw_1\wedge d\bar w_1$ so 
$$d(d\rho+i\sigma_1-\bar w_1dw_1)=0$$
and since $d\rho+i\sigma_1$ is of type $(1,0)$ we have a holomorphic function $w_2$ such that $idw_2/2=d\rho+i\sigma_1-\bar w_1dw_1$. It follows that $w_1,w_2$ are complex coordinates and moreover 
$$2d\rho=d\vert w_1\vert^2-\im w_2$$
and choosing the constant ambiguity of $w_2$ appropriately $\rho=(\vert w_1\vert^2-\im w_2)/2$. But $0<\rho<\infty$ so the complex structure is the {\it exterior} of the Siegel domain. Up to a scale we can write the metric in these coordinates as 
$$(\bar w_1dw_1-idw_2/2)\wedge(w_1d\bar w_1+id\bar w_2/2)+(\vert w_1\vert^2-\im w_2)dw_1\wedge d\bar w_1.$$

\begin{rmks} 

\noindent 1. The authors of \cite{ACDM} calculate a one-parameter family of metrics by adding a constant $c$ to  the moment map, or equivalently changing the lift of the $\R$-action. The scalar-flat metric is then 
$$g=\frac{1}{4}\left(\frac{\rho+2c}{\rho+c}d\rho^2 +\frac{\rho+c}{\rho+2c} \sigma_1^2+2(\rho+2c)(\sigma_2^2+\sigma_3^2)\right).$$
 The complex structure is the same but with $du+i\sigma_1=-\bar w_1dw_1+idw_2/2$ where $u=(\rho+c)+c\log(\rho+c)$, and the K\"ahler potential is then $f(\rho)=(\rho+c)^2+4c(\rho+c)+2c^2\log(\rho+c)$.

 \noindent 2. From (\ref{scal}) it is clear that a curve orthogonal to the $H$-orbits has finite length as $\rho\rightarrow 0$ and so the metric is not complete. The function $\rho$ is invariant under the circle action and so is defined by a  function  on the hyperk\"ahler manifold -- in fact $\vert w\vert^2$. 
\end{rmks}

\subsection{Higher dimensions}\label{42}

The hyperk\"ahler picture in higher dimensions is a straightforward generalization of the above. We take  $M=\C^m\oplus j\C^m$ and $$\omega_1=\frac{i}{2}\sum_{\alpha}(dz_{\alpha} \wedge d\bar z_{\alpha}+dw_{\alpha} \wedge d\bar w_{\alpha}), \qquad \omega_2+i\omega_3=\sum_{\alpha}dz_{\alpha}\wedge dw_{\alpha}$$
and the action  $(z,w)\mapsto (z,e^{i\theta}w)$ on $(z,w)\in \C^m\times \C^m$. Now we shall only identify the manifold of holonomy $U^*(2m)$ corresponding to this as a complex manifold and this requires that we revisit the twistor approach, as discussed in Sections 4.2 and 4.3 of \cite{Hit}. The advantage here is that we can pass directly from the hyperk\"ahler manifold with its complex structure $I$ to the $U^*(2m)$-manifold. 

The twistor space $Z$ of a hyperk\"ahler manifold $M$ is a holomorphic fibre bundle over $\CP^1$, and the circle action that we are considering induces an action which respects the fibration but acts as a rotation on the sphere $\CP^1$ leaving fixed two points $0,\infty$. The fibres over these two points are biholomorphically equivalent to the manifold $M$ with complex structures $\pm I$. Since the curvature form $F=\omega_1+dd_1^c\mu$ is of type $(1,1)$ with respect to all complex structures $I,J,K$ it defines a holomorphic principal $\C^*$-bundle $P^c$ over $Z$.  If the circle action on $M$ preserving $I$ extends to a holomorphic $\C^*$-action, lifting to the principal bundle, then in \cite{Hit} it is argued that the twistor space of the $U^*(2m)$-manifold is the quotient $P^c/\C^*$ and the two fibres over $0$ and $\infty$ give the two components of the divisor of $K^{1/(m+1)}$ which defines the complex structure $I$. Thus to determine the complex structure we just need to take the quotient of $P^c$ restricted to the fibre over $0$ by the $\C^*$-action, forgetting the intermediate quaternionic K\"ahler manifold. 

It is well-known, however, that to get a Hausdorff manifold as a holomorphic quotient we must usually choose an open set of stable points for the action. In this case it means finding the points in the trivial $\C^*$-bundle  $\C^m\times\C^m\backslash\{0\}\times \C^*$ which are equivalent under the $\C^*$-action to points in the unit circle bundle for the Hermitian metric   $\exp(\Vert z\Vert^2-\Vert w\Vert^2)/2$. 

So if $(z,w,u)\in \C^m\times\C^m\backslash\{0\}\times \C^*$ lies in the circle bundle then $\vert u\vert^2=\exp((\Vert w\Vert^2-\Vert z\Vert^2)/2)$ and transforming by $t>0$ in $\C^*$ we have 
$(z,tw,u)$. Setting $\tilde w=tw$ this means  
$$\vert u\vert^2=\exp((t^{-2}\Vert \tilde w\Vert^2-\Vert z\Vert^2)/2)$$
so that $\vert u\vert^2>\exp(-\Vert z\Vert^2/2)$ which is a proper open subset of $\C^m\times\C^m\backslash\{0\}\times \C^*$. 

Writing $u=e^{iu_2}$ and $z=2u_1$, the inequality gives $\Vert u_1\Vert^2> \im u_2$ which defines the exterior $S$ of an $(m+1)$-dimensional Siegel domain.  Taking the quotient of $\C^m\times\C^m\backslash\{0\}\times \C^*$ by $\C^*$ therefore gives the quotient of the  product 
$$\CP^{m-1}\times S$$
by the integer action $u_2\mapsto u_2+2\pi n$.  

\begin{rmks}

\noindent 1. The physicist's c-map constructs  a quaternionic K\"ahler manifold of dimension $4m$ from a projective special K\"ahler manifold of complex dimension $(m-1)$ and it is a fibration over the this space \cite{Hit1}. The authors of \cite{ACDM}  give some explicit formulas for the c-map construction and we have used the simplest one in Section \ref{41}. It seems likely, given the product structure in the flat case above, that in general the associated $U^*(2m)$-manifold fibres holomorphically over the projective special K\"ahler manifold which is $\CP^{m-1}$ in our case (see \cite{Cort}).

\noindent 2. The four-dimensional calculation suggests that the higher-dimensional transforms of flat space are not geodesically complete. However, the various types of asymptotically flat hyperk\"ahler manifolds that exist in the literature may well generate complete $U^*(2m)$-manifolds. Certainly the 4-dimensional Calabi metric, otherwise known as the Eguchi-Hanson metric,  is asymptotically locally Euclidean and transforms, as we saw,  to $D\times \CP^1$ which is complete.

\noindent 3. An isometry group $G$ of the hyperk\"ahler manifold which commutes with the circle action will generate an action of a central extension   on the quaternionic K\"ahler manifold and correspondingly of the $U^*(2m)$-manifold. Using the flat starting point this should give scope for the application of Joyce's quotient construction \cite{Joy1} to generate more examples. 
\end{rmks}

 \end{document}